\documentclass{amsart}
\usepackage{amsfonts}
\usepackage{amsmath}
\usepackage{amssymb}
\usepackage{amsthm}

\newtheorem{theorem}{Theorem}[section]
\newtheorem{prop}{Proposition}[section]
\newtheorem{lemma}{Lemma}[section]
\newtheorem{rem}{Remark}[section]
\newtheorem{exmp}{Example}[section]

\begin{document}
\author{Mark Pankov}
\title{Isometric embeddings of half-cube graphs in half-spin Grassmannians}
\subjclass[2000]{}
\address{Department of Mathematics and Computer Science, University of Warmia and Mazury,
S{\l}oneczna 54, Olsztyn, Poland}
\email{pankov@matman.uwm.edu.pl, markpankovmath@gmail.com}

\begin{abstract}
Let $\Pi$ be a polar space of type $\textsf{D}_{n}$.
Denote by ${\mathcal G}_{\delta}(\Pi)$, $\delta\in \{+,-\}$ the associated half-spin Grassmannians
and write $\Gamma_{\delta}(\Pi)$ for the corresponding half-spin Grassmann graphs.
In the case when $n\ge 4$ is even, the apartments of ${\mathcal G}_{\delta}(\Pi)$
will be characterized as the images of isometric embeddings of 
the half-cube graph $\frac{1}{2}H_n$ in $\Gamma_{\delta}(\Pi)$.
As an application,
we describe all isometric embeddings of $\Gamma_{\delta}(\Pi)$
in the half-spin Grassmann graphs associated to a polar space of type $\textsf{D}_{n'}$
under the assumption that $n\ge 6$ is even.
\end{abstract}

\maketitle
\section{Introduction}

In the present paper we continue to discuss 
the problem of metric characterization of apartments in building Grassmannians \cite{Pankov2,Pankov3}.
This problem is connected to the results obtained in \cite{CKS,Kasikova1,Kasikova2,Pankov4}.

By \cite{Tits}, a {\it building} is a simplicial complex $\Delta$ containing a family of subcomplexes
called {\it apartments} and characterized by some properties.
All apartments are isomorphic to a certain Coxeter complex,
i.e. the simplicial complex associated to a Coxeter system,
 which defines the type of the building.
We suppose that our building is {\it spherical}, i.e. 
the associated Coxeter system is finite.
Maximal simplices of $\Delta$ are said to be {\it chambers}.
They have the same finite cardinality $n$ called the {\it rank} of $\Delta$.
%Two chambers are {\it adjacent} if their intersection consists of $n-1$ vertices.
The vertex set of $\Delta$ can be labeled by the nodes of the diagram corresponding to 
the associated Coxeter system (such labeling is unique up to a permutation on the set of nodes).
Every set consisting of all vertices labeled by the same node 
is said to be a {\it Grassmannian}.
The intersection of every chamber and every Grassmannian is a single vertex.
Thus the vertex set of $\Delta$ is decomposed  in precisely $n$ distinct Grassmannians.
The intersections of apartments with one of the building Grassmannians are called 
{\it apartments in this Grassmannian}.

Let ${\mathcal G}$ be one of the building Grassmannians.
We say that $a,b\in {\mathcal G}$ are {\it adjacent} if
there is a simplex $P$ such that $P\cup\{a\}$ and $P\cup\{b\}$ are chambers.
Consider the associated {\it Grassmann graph} $\Gamma$, 
i.e. the graph whose vertex set is ${\mathcal G}$ and whose edges are pairs of adjacent vertices.
Suppose that ${\mathcal A}$ is an apartment of ${\mathcal G}$ and 
denote by $\Gamma_{\mathcal A}$ the restriction of the graph $\Gamma$ to ${\mathcal A}$.
We want to distinguish all cases such that the image of every isometric embedding of $\Gamma_{\mathcal A}$
in $\Gamma$ is an apartment of ${\mathcal G}$. 
Note that in some cases this does not hold \cite{CS,Pankov2}.

By \cite{Tits}, there are the following seven types of {\it irreducible thick spherical buildings} of rank $\ge 3$:
three classical types
${\textsf A}_{n}, {\textsf B}_{n}={\textsf C}_{n}, {\textsf D}_{n}$
and four exceptional types ${\textsf F}_{4}, {\textsf E}_{i},i=6,7,8$.

Every building of type ${\textsf A}_{n-1}$ is the flag complex of a certain $n$-dimensional vector space $V$ 
(over a division ring). 
The associated Grassmannians are ${\mathcal G}_{k}(V)$, where $k\in\{1,\dots,n-1\}$.
Each ${\mathcal G}_{k}(V)$ consists of all $k$-dimensional subspaces of $V$
and two elements of ${\mathcal G}_{k}(V)$ are adjacent if their intersection is
$(k-1)$-dimensional. The corresponding Grassmann graph is denoted by $\Gamma_{k}(V)$.
The case when $k=1,n-1$ is trivial --- any two distinct vertices of $\Gamma_{k}(V)$ are adjacent.
Every apartment of ${\mathcal G}_{k}(V)$  is defined by a certain base of $V$:
it consists of all $k$-dimensional subspaces spanned by subsets of this base.
All apartments of ${\mathcal G}_{k}(V)$ are the images of isometric embeddings of 
the Johnson graph $J(n,k)$ in $\Gamma_{k}(V)$.
However, the image of every isometric embedding of $J(n,k)$ in $\Gamma_{k}(V)$
is an apartment if and only if $n=2k$ \cite{Pankov2}.

All buildings of types $\textsf{C}_{n}$ and $\textsf{D}_{n}$ are defined by polar spaces.
Every building of type $\textsf{C}_{n}$ is the flag complex formed by singular subspaces of a rank $n$ polar space.
The associated Grassmannians are the {\it polar Grassmannians}
consisting  of singular subspaces of the same dimension.
The polar Grassmannian whose elements are maximal singular subspaces is called the {\it dual polar space}
and the associated Grassmann graph is known as the {\it dual polar graph}.
By \cite{Pankov3}, the apartments in the dual polar space can be characterized as 
the images of isometric embeddings of the $n$-dimensional hypercube graph $H_{n}$ in the dual polar graph.

Every building of type $\textsf{D}_{n}$ can be obtained from a polar space of type $\textsf{D}_{n}$.
This construction is known as the {\it oriflamme complex}.
The Grassmannians of this building are some of the polar Grassmannians and 
so-called {\it half-spin Grassmannians}.
The Grassmann graphs associated  to the half-spin Grassmannians are called 
the {\it half-spin Grassmann graphs}.
In this paper we show that the apartments of the half-spin Grassmannians 
can be characterized as the images of isometric embeddings of 
the half-cube graph $\frac{1}{2}H_n$ in the half-spin Grassmann graphs if $n\ge 4$ is even
(Theorem \ref{theorem-apart-halfspin}).
In the case when $n$ is odd, we conjecture the existence 
of isometric embeddings of $\frac{1}{2}H_n$ in the half-spin Grassmann graphs
whose images are not apartments (Section 6).
As an application of the main result,
we describe all isometric embeddings of 
the half-spin Grassmann graphs of a polar space of type $\textsf{D}_{n}$,
where $n\ge 6$ is even, 
in the half-spin Grassmann graphs associated to a polar space of type $\textsf{D}_{n'}$
(Theorem \ref{theorem-half-spin-emb}).

Note that in \cite{CKS} apartments in Grassmannians of finite-dimensional vector spaces, 
dual polar spaces and half-spin Grassmannians
were characterized in terms of independent subsets in the associated partial linear spaces.
Some more general results can be found in \cite{Kasikova1,Kasikova2}.

\section{Basic definitions and constructions}

\subsection{Graphs}
The {\it distance} between two vertices in a connected graph is
the smallest number $i$ such that there is a path of length $i$
(a path consisting of $i$ edges) connecting  these vertices.
A path connecting two vertices
is said to be a {\it geodesic} if the number of edges in this path
is equal to the distance between the vertices.
The maximum of all distances between vertices in a graph
is called the {\it diameter} of the graph.
In the case when the diameter is finite,
two vertices are said to be {\it opposite}
if the distance between them is equal to the diameter.

A subset in the vertex set of a graph is called a {\it clique} if any two distinct elements of this subset
are adjacent vertices. 
Maximal cliques exist and every clique is contained in
a certain maximal clique. For finite graph this is trivial
and we use Zorn lemma for infinite graphs.

An injective mapping of the vertex set of a graph $\Gamma$ to the vertex set of a graph $\Gamma'$ 
is an {\it embedding} of $\Gamma$ in $\Gamma'$ if 
it transfers adjacent vertices to adjacent vertices and non-adjacent vertices to non-adjacent vertices.
Every surjective embedding is an isomorphism between the graphs.
An embedding is said to be {\it isometric} if it preserves
the distance between vertices.

\subsection{Hypercube and half-cube graphs}
The vertex set of the $n$-dimensional {\it hypercube graph} $H_{n}$ is formed by all
sequences $(a_{1},\dots,a_{n})$, where each $a_{i}$ is equal to $0$ or $1$.
For any two such sequences $a=(a_{1},\dots,a_{n})$ and $b=(b_{1},\dots,b_{n})$
we define the distance 
$$d(a,b):=|a_{1}-b_{1}|+\dots+|a_{n}-b_{n}|.$$
The sequences $a$ and $b$ are adjacent vertices of  $H_{n}$ if $d(a,b)=1$.
The graph $H_{n}$ is connected, $d$ is the distance on $H_{n}$
and the diameter of $H_{n}$ is equal to $n$.

\begin{lemma}\label{lemma-cube}
For  every vertex of $H_{n}$ there is unique vertex of $H_{n}$ opposite to this vertex.
If $v$ and $w$ are opposite vertices of $H_{n}$ then every vertex of $H_{n}$
belongs to a geodesic connecting $v$ and $w$.
\end{lemma}

\begin{proof}
 Easy verification.
\end{proof}

\begin{rem}\label{rem-cube}{\rm
Consider the $(2n)$-element subset $J=\{\pm 1,\dots,\pm n\}$.
A subset $I\subset J$ is said to be {\it singular} 
if for every $i\in I$ we have $-i\not\in I$.
Every maximal singular subset consists of $n$ elements.
The hypercube graph $H_{n}$ can be defined 
as the graph whose vertex set is formed by all maximal singular subsets.
Two such subsets are adjacent vertices of $H_{n}$ if their intersection consists of $n-1$ elements.
}\end{rem}

The vertex set of $H_{n}$ can be uniquely decomposed in  the sum of two disjoint subsets 
$A_{+}$ and $A_{-}$ such that the distance between any two elements of $A_{\delta}$, $\delta\in \{+,-\}$ 
is even and the distance between every element of $A_{+}$ and every element of $A_{-}$ is odd.
Consider the graph  $H_{\delta}$, $\delta\in \{+,-\}$ whose vertex set is $A_{\delta}$
and two vertices of $H_{\delta}$ are adjacent if the distance between them 
in $H_{n}$ is equal to $2$. 
The graphs $H_{+}$ and $H_{-}$ are isomorphic.
Each of these graphs is called the $n$-dimensional {\it half-cube graph} and 
denoted by $\frac{1}{2}H_{n}$.

The graph $\frac{1}{2}H_{n}$ is connected.
If $v$ and $u$  are vertices of $\frac{1}{2}H_{n}$
then the distance between them is equal to $\frac{d(v,u)}{2}$.
The diameter of  $\frac{1}{2}H_{n}$ is equal to $\frac{n}{2}$ or $\frac{n-1}{2}$
if $n$ is even or odd, respectively.

There is the following analogue of Lemma \ref{lemma-cube}.

\begin{lemma}\label{lemma-half-cube-even}
Suppose that $n$ is even.
For  every vertex of $\frac{1}{2}H_{n}$ there is unique vertex of 
$\frac{1}{2}H_{n}$ opposite to this vertex.
If $v$ and $w$ are opposite vertices of $\frac{1}{2}H_{n}$ then every vertex of 
$\frac{1}{2}H_{n}$
belongs to a geodesic connecting $v$ and $w$.
\end{lemma}

\begin{proof}
 Easy verification.
\end{proof}

If $n$ is odd then the above statement fails and
for every vertex of $\frac{1}{2}H_{n}$ there are precisely $n$ vertices 
of $\frac{1}{2}H_{n}$ opposite to this vertex.

\subsection{Partial linear spaces}
Consider a non-empty set $P$ whose elements are called {\it points}
and a family ${\mathcal L}$ of proper subsets of $P$ called {\it lines}.
We say that two or more points are {\it collinear} if there is a line containing all of them.
Suppose that the pair $\Pi=(P,{\mathcal L})$ is a {\it partial linear space},
i.e. the following axioms hold:
\begin{enumerate}
\item[$\bullet$]
every line contains at least two points and every point belongs to a line,
\item[$\bullet$]
for any distinct collinear points $p,q\in P$ there is precisely one line $pq$ containing them.
\end{enumerate}
A subset $S\subset P$ is a {\it subspace} of $\Pi$
if for any distinct collinear points $p,q\in S$
the line $pq$ is contained in $S$.
A subspace is said to be {\it singular} if any two distinct points of the subspace are collinear
(by the definition, the empty set and a single point are singular subspaces).
Using Zorn lemma, we establish the existence of maximal singular subspaces and the fact that 
every singular subspace is contained in a certain maximal singular subspace.

For every subset $X\subset P$ the minimal subspace containing $X$,
i.e. the intersection of all subspaces containing $X$, is called {\it spanned} by $X$
and denoted by $\langle X\rangle$.
We say that $X$ is {\it independent} if the subspace $\langle X\rangle$ 
cannot be spanned by a proper subset of $X$.

Let $S$ be a subspace of $\Pi$ (possible $S = P$).
An independent subset $X\subset S$ is a {\it base} of $S$ if $\langle X\rangle=S$.
The {\it dimension} of $S$ is the smallest cardinality $\alpha$
such that $S$ has a base of cardinality $\alpha+1$.
The dimension of the empty set and a single point is equal to $-1$ and $0$ (respectively),
lines are $1$-dimensional subspaces.
Every $2$-dimensional singular subspace will be called a {\it plane}.

Two partial linear spaces $\Pi=(P,{\mathcal L})$ and $\Pi'=(P',{\mathcal L}')$
are {\it isomorphic} if there is a bijection $f:P\to P'$ satisfying 
$f({\mathcal L})={\mathcal L}'$. Every such bijection is called a {\it collineation}
of $\Pi$ to $\Pi'$.

\subsection{Polar spaces}
By \cite{BuekenhoutCohen-book,Pankov1,Shult-book,Ueberberg}, 
a {\it polar space} is a partial linear space $\Pi=(P,{\mathcal L})$
satisfying the following axioms:
\begin{enumerate}
\item[$\bullet$] every line contains at least three points,
\item[$\bullet$] there is no point collinear to all points,
\item[$\bullet$] if $p\in P$ and $L\in {\mathcal L}$
then $p$ is collinear to one or all points of the line $L$,
\item[$\bullet$] any flag formed by singular subspaces is finite.
\end{enumerate}
If our polar space contains a singular subspace whose dimension is not less than $2$
then all maximal singular subspaces of $\Pi$ are projective spaces
of the same dimension $n\ge 2$ and the number $n+1$ is called the {\it rank} of $\Pi$.

The collinearity relation of $\Pi$ will be denoted by $\perp$.
For points $p,q\in P$ we write $p\perp q$ if $p$ is collinear to $q$
and $p\not\perp q$ otherwise.
Moreover, if $X,Y\subset P$ then $X\perp Y$ means that every point of $X$
is collinear to all points of $Y$.
If $X\perp X$ then the subspace $\langle X\rangle$ is singular.

For every polar space of rank $n$  one of the following possibilities is realized:
\begin{enumerate}
\item[$\bullet$] type ${\mathsf C}_{n}$ ---
every $(n-2)$-dimensional singular subspace is contained in at least three maximal singular subspaces,
\item[$\bullet$] type ${\mathsf D}_{n}$ ---
every $(n-2)$-dimensional singular subspace is contained in precisely two maximal singular subspaces.
\end{enumerate}
All polar spaces of rank $\ge 3$ were described by J. Tits \cite{Tits}.
We will focus our attention on polar spaces of type ${\mathsf D}_{n}$, $n\ge 4$.

\begin{exmp}\label{exmp-polar-D}{\rm
Let $V$ be a  $(2n)$-dimensional vector space over a field.
Suppose that the characteristic of this field is not equal to $2$ 
and consider a non-degenerate symmetric bilinear form on $V$ 
such that maximal totally isotropic subspaces are $n$-dimensional.
The associated polar space whose points are $1$-dimensional isotropic subspaces
and whose lines are defined by $2$-dimensional totally isotropic subspaces
(the line corresponding to a $2$-dimensional totally isotropic subspace $S$
consists of all $1$-dimensional subspaces of $S$)
is of type ${\mathsf D}_{n}$.
In the case when the characteristic of the field is equal to $2$,
we consider a non-defect quadratic form on $V$
such that maximal singular subspaces are $n$-dimensional.  
The associated polar space 
(the points are $1$-dimensional singular subspaces
and the lines are defined by $2$-dimensional singular subspaces)
is also of type ${\mathsf D}_{n}$.
}\end{exmp}

It follows from Tits's description of polar spaces 
that every polar space of type ${\mathsf D}_{n}$, $n\ge 4$ is isomorphic to 
one of the polar spaces presented in Example \ref{exmp-polar-D}.

\subsection{Polar Grassmannians}
Let $\Pi=(P,{\mathcal L})$ be a polar space of rank $n$.
For every $k\in\{0,1,\dots,n-1\}$ we denote by ${\mathcal G}_{k}(\Pi)$
the polar Grassmannian consisting of all $k$-dimensional singular subspaces of $\Pi$;
in particular, ${\mathcal G}_{n-1}(\Pi)$ is formed by maximal singular subspaces.

The {\it dual polar graph} $\Gamma_{n-1}(\Pi)$ is the graph
whose vertex set is ${\mathcal G}_{n-1}(\Pi)$ and
two distinct elements of ${\mathcal G}_{n-1}(\Pi)$ are adjacent vertices of $\Gamma_{n-1}(\Pi)$
if their intersection is an $(n-2)$-dimensional singular subspace.
The graph $\Gamma_{n-1}(\Pi)$ is connected and 
the distance $d(S,U)$ between 
$S,U\in {\mathcal G}_{n-1}(\Pi)$ is equal to
$$n-1-\dim(S\cap U).$$
The diameter of $\Gamma_{n-1}(\Pi)$ is equal to $n$ and two vertices of $\Gamma_{n-1}(\Pi)$
are opposite if and only if they are disjoint singular subspaces.

Let $S$ and $U$ be singular subspaces of $\Pi$ such that $S\subset U$.
If $$\dim S<k<\dim U$$ then 
we write $[S,U]_{k}$ for the set of all $X\in {\mathcal G}_{k}(\Pi)$
satisfying $S\subset X\subset U$.
In the case when $S=\emptyset$, this set will be denoted by $\langle U]_{k}$.

Let $M$ be an $m$-dimensional singular subspace of $\Pi$.
For every natural $k>m$ we denote by $[M\rangle_{k}$ the set of all elements of 
${\mathcal G}_{k}(\Pi)$ containing $M$.
Now suppose that $m<n-2$. 
For every $N\in {\mathcal G}_{m+2}(\Pi)$ containing $M$
the set $[M,N]_{m+1}$ is called a {\it line} of $[M\rangle_{m+1}$.
The set  $[M\rangle_{m+1}$ together with the family of all such lines 
is a polar space of rank $n-m-1$ \cite[Lemma 4.4]{Pankov1}.
This polar space will be denoted by $\Pi_{M}$.
Every maximal singular subspace of $\Pi_{M}$ is of type $[M,U]_{m+1}$,
where $U$ is a maximal singular subspace of $\Pi$ containing $M$.
Thus $[M\rangle_{n-1}$ can be naturally identified with ${\mathcal G}_{n-m-2}(\Pi_{M})$ 
and the dual polar graph $\Gamma_{n-m-2}(\Pi_{M})$ is the restriction of 
the graph $\Gamma_{n-1}(\Pi)$ to $[M\rangle_{n-1}$.
In the case when $\Pi$ is a polar space of type $\textsf{D}_{n}$,
the polar space $\Pi_{M}$ is of type $\textsf{D}_{n-m-1}$.

\subsection{Half-spin Grassmannians}
Let $\Pi=(P,{\mathcal L})$ be a polar space of type $\textsf{D}_{n}$.
Then ${\mathcal G}_{n-1}(\Pi)$ can be uniquely decomposed in the sum of two
disjoint subsets, we denote them by ${\mathcal G}_{+}(\Pi)$ and ${\mathcal G}_{-}(\Pi)$, such that
the distance between any two elements of ${\mathcal G}_{\delta}(\Pi)$, $\delta\in \{+,-\}$ 
in the dual polar graph $\Gamma_{n-1}(\Pi)$
is even  and the same distance between any $S\in {\mathcal G}_{\delta}(\Pi)$ and $U\in {\mathcal G}_{-\delta}(\Pi)$ 
is odd
(we write $-\delta$ for the complement of $\delta$ in the set $\{+,-\}$).
These subsets are called the {\it half-spin Grassmannians} of $\Pi$.

The {\it half-spin Grassmann graph} $\Gamma_{\delta}(\Pi)$, $\delta\in \{+,-\}$
is the graph whose vertex set is ${\mathcal G}_{\delta}(\Pi)$
and two elements of ${\mathcal G}_{\delta}(\Pi)$ are adjacent vertices of $\Gamma_{\delta}(\Pi)$
if their intersection is an $(n-3)$-dimensional singular subspace, 
i.e. the distance between the corresponding vertices of $\Gamma_{n-1}(\Pi)$ is equal to $2$.

This graph is connected and the distance $d_{\delta}(S,U)$ between 
$S,U\in {\mathcal G}_{\delta}(\Pi)$ in $\Gamma_{\delta}(\Pi)$
is equal to $\frac{d(S,U)}{2}$.
If $n$ is even then the diameter of $\Gamma_{\delta}(\Pi)$ 
is equal to $\frac{n}{2}$ and two vertices are opposite
if and only if they are disjoint singular subspaces. 
In the case when $n$ is odd,
the diameter is equal to $\frac{n-1}{2}$ and 
two vertices are opposite if and only if their intersection is a single point.

\begin{rem}\label{rem-halfspin}{\rm
If a polar space is defined by a non-degenerate symmetric bilinear form $\Omega$
(see Example \ref{exmp-polar-D}) 
then their maximal singular subspaces can be identified with 
the maximal totally isotropic subspaces of $\Omega$
and the half-spin Grassmannians are the orbits of the action of the orthogonal group ${\rm O}^{+}(\Omega)$
on the set of all maximal totally isotropic subspaces.
Every element of ${\rm O}(\Omega)\setminus {\rm O}^{+}(\Omega)$ transfers 
one half-spin Grassmannian to the other and it induces an isomorphism between 
the half-spin Grassmann graphs. 
The same holds for the polar spaces defined by quadratic forms.
}\end{rem}

The latter remark shows that the graphs $\Gamma_{+}(\Pi)$ and  $\Gamma_{-}(\Pi)$
are isomorphic.
For $n=2,3$ any two distinct vertices of $\Gamma_{\delta}(\Pi)$ are adjacent 
and we will always suppose that $n\ge 4$.

Let $M$ be an $m$-dimensional singular subspace of $\Pi$ and $m< n-2$.
Denote by $[M\rangle_{\delta}$ 
the set of all elements of ${\mathcal G}_{\delta}(\Pi)$ containing $M$.
In the case when $m=n-3$, this set is called a {\it line} of ${\mathcal G}_{\delta}(\Pi)$.
Any two distinct elements of such a line are adjacent vertices of $\Gamma_{\delta}(\Pi)$.
If $m<n-3$ then $\Pi_{M}$ is a polar space of type $\textsf{D}_{n-m-1}$
whose half-spin Grassmannian ${\mathcal G}_{\delta}(\Pi_{M})$
can be naturally identified with $[M\rangle_{\delta}$ and 
the half-spin Grassmann graph $\Gamma_{\delta}(\Pi_{M})$ is
the restriction of the graph $\Gamma_{\delta}(\Pi)$ to $[M\rangle_{\delta}$.

\section{Main result}
Let $\Pi=(P,{\mathcal L})$ be a polar space of rank $n$.
A set ${\mathcal F}$ consisting of $2n$ distinct points $p_{1},\dots,p_{2n}\in P$ 
is called a {\it frame} of $\Pi$
if for every $i\in\{1\,\dots,2n\}$ there exists unique $\sigma(i)\in\{1\,\dots,2n\}$
such that $p_{i}\not\perp p_{\sigma(i)}$.
This is an  independent subset  and any $k$ distinct mutually collinear points of ${\mathcal F}$
span a $(k-1)$-dimensional singular subspace.
Denote by ${\mathcal A}$ the set consisting of all maximal singular 
subspaces spanned by subsets of ${\mathcal F}$, i.e.
all subspaces of type $$\langle p_{i_{1}},\dots, p_{i_{n}}\rangle$$ such that
$$\{i_{1},\dots, i_{n}\}\cap \{\sigma(i_{1}),\dots, \sigma(i_{n})\}=\emptyset.$$
This set is  the {\it apartment} of ${\mathcal G}_{n-1}(\Pi)$ associated to the frame ${\mathcal F}$.
It follows from Remark \ref{rem-cube} that
${\mathcal A}$ is the image of an isometric embedding of
the $n$-dimensional hypercube graph $H_{n}$ in 
the dual polar graph $\Gamma_{n-1}(\Pi)$.
If $M$ is an $(n-m-1)$-dimensional singular subspace of $\Pi$
then $\Pi_{M}$ is a polar space of rank $m$ 
and every apartment of ${\mathcal G}_{m-1}(\Pi_{M})=[M\rangle_{n-1}$
is the image of an isometric embedding of $H_{m}$ in $\Gamma_{n-1}(\Pi)$.

\begin{theorem}[\cite{Pankov3}]\label{theorem-apart-dual}
The image of every isometric embedding of 
the $m$-dimensional hypercube graph $H_{m}$, $m\le n$
in the dual polar graph $\Gamma_{n-1}(\Pi)$ is an apartment of 
${\mathcal G}_{m-1}(\Pi_{M})$, where
$M$ is an $(n-m-1)$-dimensional singular subspace of $\Pi$.
In particular, the image of every isometric embedding of $H_{n}$
in $\Gamma_{n-1}(\Pi)$ is an apartment of ${\mathcal G}_{n-1}(\Pi)$.
\end{theorem}

Now suppose that $\Pi$ is a polar space of type $\textsf{D}_{n}$
and $\delta\in\{+,-\}$.
The intersection 
$${\mathcal A}_{\delta}:={\mathcal A}\cap {\mathcal G}_{\delta}(\Pi)$$
is the {\it apartment} of the half-spin Grassmannian ${\mathcal G}_{\delta}(\Pi)$ 
associated to the frame ${\mathcal F}$.
This is the image of an isometric embedding of 
the half-cube graph $\frac{1}{2}H_{n}$ in the half-spin Grassmann graph $\Gamma_{\delta}(\Pi)$
(this easy follows from the fact that 
${\mathcal A}$ is the image of an isometric embedding of $H_{n}$ in $\Gamma_{n-1}(\Pi)$).
If $M$ is an $(n-m-1)$-dimensional singular subspace of $\Pi$
then $\Pi_{M}$ is a polar space of type $\textsf{D}_{m}$
and every apartment of ${\mathcal G}_{\delta}(\Pi_{M})=[M\rangle_{\delta}$
is the image of an isometric embedding of $\frac{1}{2}H_{m}$ in $\Gamma_{\delta}(\Pi)$.

\begin{exmp}\label{exmp3-1}{\rm
Suppose that $n=4$. Then ${\mathcal G}_{\delta}(\Pi)$ together with 
the family of all lines is a polar space of type $\textsf{D}_{4}$ \cite[Proposition 4.23]{Pankov1}.
This polar space will be denoted by $\Pi_{\delta}$.
Two distinct points of $\Pi_{\delta}$ are collinear if and only if they are
adjacent vertices of $\Gamma_{\delta}(\Pi)$.
The graph $\frac{1}{2}H_{4}$ consists of $8$ vertices and 
for every vertex there is precisely one vertex non-adjacent to it.
This means that the frames of $\Pi_{\delta}$ can be characterized as
the images of  embeddings of $\frac{1}{2}H_{4}$ in $\Gamma_{\delta}(\Pi)$
(every such embedding is isometric, since the both graphs are of diameter $2$).
By \cite[Corollary 4.4]{Pankov1} (see also \cite{CKS}),
the family of all apartments of ${\mathcal G}_{\delta}(\Pi)$ coincides with the family of all frames of $\Pi_{\delta}$.
So, the image of every embedding of $\frac{1}{2}H_{4}$ in $\Gamma_{\delta}(\Pi)$
is an apartment of ${\mathcal G}_{\delta}(\Pi)$.
}\end{exmp}

The main result of the present paper is the following.

\begin{theorem}\label{theorem-apart-halfspin}
Suppose that $\Pi$ is a polar space of type {\rm $\textsf{D}_{n}$}.
If $m$ is an even integer satisfying $4\le m\le n$
then the image of every isometric embedding of
the $m$-dimensional half-cube graph $\frac{1}{2}H_{m}$
in the half-spin Grassmann graph $\Gamma_{\delta}(\Pi)$, $\delta\in\{+,-\}$
is an apartment of ${\mathcal G}_{\delta}(\Pi_{M})$, where
$M$ is an $(n-m-1)$-dimensional singular subspace of $\Pi$.
In particular,  if $n$ is even then the image of every isometric embedding of $\frac{1}{2}H_{n}$
in $\Gamma_{\delta}(\Pi)$ is an apartment of ${\mathcal G}_{\delta}(\Pi)$.
\end{theorem}

In Section 7, we apply Theorem \ref{theorem-apart-halfspin}
to isometric embeddings of half-spin Grassmann graphs.

\section{Cliques}

\subsection{Maximal cliques of half-spin Grassmannian graphs}
Let $\Pi=(P,{\mathcal L})$ be a polar space of type $\textsf{D}_{n}$
and $\delta\in\{+,-\}$. 
By \cite[Subsection 4.5.2]{Pankov1},
there are precisely the following two types of maximal cliques of $\Gamma_{\delta}(\Pi)$:
\begin{enumerate}
\item[$\bullet$] The {\it star} $[M\rangle_{\delta}$, $M\in {\mathcal G}_{n-4}(\Pi)$
(if two distinct elements of ${\mathcal G}_{\delta}(\Pi)$ contain $M$
then their intersection is $(n-3)$-dimensional and they are adjacent vertices of  $\Gamma_{\delta}(\Pi)$).
The star $[M\rangle_{\delta}$ together with all lines of ${\mathcal G}_{\delta}(\Pi)$
contained in $[M\rangle_{\delta}$ is a $3$-dimensional projective space.
\item[$\bullet$] The {\it special subspace} $[U]_{\delta}$, $U\in {\mathcal G}_{-\delta}(\Pi)$ 
formed by all elements of ${\mathcal G}_{\delta}(\Pi)$ intersecting $U$ in 
$(n-2)$-dimensional singular subspaces, i.e.
all vertices of $\Gamma_{n-1}(\Pi)$ adjacent to $U$.
The  special subspace $[U]_{\delta}$ together with all lines of ${\mathcal G}_{\delta}(\Pi)$
contained in $[U]_{\delta}$ is an $(n-1)$-dimensional projective space.
\end{enumerate}
Recall the following facts concerning the intersection of two distinct maximal cliques:
\begin{enumerate}
\item[$\bullet$]
The intersection of two distinct stars $[M\rangle_{\delta}$ and $[M'\rangle_{\delta}$
is empty or a single vertex or a line. The second possibility is not realized  if $n=4$.
The intersection is a line if and only if $M$ and $M'$ span an $(n-3)$-dimensional singular subspace.
\item[$\bullet$] 
The intersection of two distinct special subspaces is empty or a line. 
The second possibility is realized if and only if the associated elements of ${\mathcal G}_{-\delta}(\Pi)$
are adjacent vertices of $\Gamma_{-\delta}(\Pi)$.
\item[$\bullet$] 
The dimension of the intersection of a star $[M\rangle_{\delta}$ and a special subspace $[U]_{\delta}$ 
is not greater than $2$. 
This intersection is a plane if and only if 
$M$ is contained in $U$.
\end{enumerate}
Therefore, the dimension of the intersection of two distinct maximal cliques of $\Gamma_{\delta}(\Pi)$
is not greater than $2$. 
If this intersection is a plane then one of the cliques is a star and the other is a special subspace.

\begin{lemma}\label{lemma4-1}
For any distinct maximal cliques ${\mathcal C}$ and ${\mathcal C}'$ in $\Gamma_{\delta}(\Pi)$
there is a sequence of maximal cliques
$${\mathcal C}={\mathcal C}_{0},{\mathcal C}_{1},\dots, {\mathcal C}_{i}={\mathcal C}'$$
such that ${\mathcal C}_{j-1}\cap {\mathcal C}_{j}$ is a plane  
for every $j\in \{1,\dots,i\}$. 
The cliques ${\mathcal C}$ and ${\mathcal C}'$ are of the same type
{\rm(}both are stars or both are special subspaces{\rm)}
 if and only if $i$ is even.
\end{lemma}

\begin{proof}
Easy verification.
\end{proof}

\subsection{Maximal cliques of half-cube graphs}
As in the previous subsection, we suppose that $\Pi=(P,{\mathcal L})$
is a polar space of type $\textsf{D}_n$.
Let $\delta\in \{+,-\}$ and let ${\mathcal A}$ be the apartment of ${\mathcal G}_{\delta}(\Pi)$
defined by a frame ${\mathcal F}$.

The restriction of $\Gamma_{\delta}(\Pi)$ to  ${\mathcal A}$ is isomorphic to 
the half-cube graph $\frac{1}{2}H_{n}$. 
Every clique of this restriction is a clique in $\Gamma_{\delta}(\Pi)$.
Therefore, this graph has precisely the following two types of maximal cliques:
\begin{enumerate}
\item[$\bullet$] The {\it star} ${\mathcal A}\cap [S\rangle_{\delta}$,
where $S\in {\mathcal G}_{n-4}(\Pi)$ is spanned by a subset of ${\mathcal F}$.
This is a base of the projective space $[S\rangle_{\delta}$,
i.e. it consists of $4$ vertices.
\item[$\bullet$] The {\it special subset} ${\mathcal A}\cap [U]_{\delta}$, 
where $U\in {\mathcal G}_{-\delta}(\Pi)$ is spanned by a subset of ${\mathcal F}$.
This is a base of the projective space $[U]_{\delta}$.
It consists of $n$ elements of ${\mathcal A}$ adjacent to $U$ in $\Gamma_{n-1}(\Pi)$.
\end{enumerate}
The intersection of two distinct stars contains at most two vertices.
This intersection is maximal  if and only if the associated elements of ${\mathcal G}_{n-4}(\Pi)$
span an element of ${\mathcal G}_{n-3}(\Pi)$.
Similarly, the intersection of two distinct special subsets contains at most two vertices 
and this intersection is maximal only in the case when the associated elements of ${\mathcal G}_{-\delta}(\Pi)$
are adjacent vertices of $\Gamma_{-\delta}(\Pi)$.
The intersection of a star and a special subset contains at most three vertices.
This intersection is maximal if and only if the associated elements of ${\mathcal G}_{n-4}(\Pi)$ 
and ${\mathcal G}_{-\delta}(\Pi)$ are incident.
Therefore, the intersection of two distinct maximal cliques contains precisely three vertices
if and only if one of these maximal cliques is a star and the other is a special subset.
There is the following analogue of Lemma \ref{lemma4-1}.

\begin{lemma}\label{lemma4-2}
For any distinct maximal cliques ${\mathcal C}$ and ${\mathcal C}'$
there is a sequence of maximal cliques
$${\mathcal C}={\mathcal C}_{0},{\mathcal C}_{1},\dots, {\mathcal C}_{i}={\mathcal C}'$$
such that $|{\mathcal C}_{j-1}\cap {\mathcal C}_{j}|=3$   
for every $j\in \{1,\dots,i\}$. 
The cliques ${\mathcal C}$ and ${\mathcal C}'$ are of the same type if and only if $i$ is even.
\end{lemma}

Let ${\mathcal A}'$ be another apartment of ${\mathcal G}_{\delta}(\Pi)$.
Suppose that $h$ is an isomorphism between the restrictions of $\Gamma_{\delta}(\Pi)$ to 
${\mathcal A}$ and ${\mathcal A}'$.
A star and a special subset are of the same cardinality only in the case when $n=4$.
This means that $h$ sends stars to stars and special subsets to special subsets if $n>4$.
If $n=4$ and there is a star whose image is a special subset
then Lemma \ref{lemma4-2} implies that all stars go to special subsets and 
all special subsets go to stars. 

In the next section we will use the following fact: 
if $h$ transfers stars to stars and special subsets to special subsets then
for every $S\in {\mathcal G}_{k}(\Pi)$, $k\in \{0,\dots,n-2\}$ 
spanned by a subset of ${\mathcal F}$
there is $S'\in {\mathcal G}_{k}(\Pi)$ spanned by a subset of the frame associated to
${\mathcal A}'$ such that 
$$h({\mathcal A}\cap [S\rangle_{\delta})={\mathcal A}'\cap [S'\rangle_{\delta}.$$

\section{Proof of Theorem \ref{theorem-apart-halfspin}}
\subsection{Reduction}
We show that the general case can be reduced to the case when $m=n$.
Let $\Pi=(P,{\mathcal L})$ be a polar space of rank $n$.

\begin{lemma}[Lemma 2 in \cite{Pankov3}]\label{lemma5-1red}
If $X_{0},\dots, X_{i}$ is a geodesic in $\Gamma_{n-1}(\Pi)$
then
$$X_{0}\cap X_{i} \subset X_{j}$$
for every $j\in\{1,\dots,i-1\}$.
\end{lemma}

Now suppose that $\Pi$ is a polar space of type $\textsf{D}_{n}$
and $\delta \in \{+,-\}$.

\begin{lemma}\label{lemma5-2red}
If $Y_{0},\dots, Y_{i}$ is a geodesic in $\Gamma_{\delta}(\Pi)$
then
$$Y_{0}\cap Y_{i} \subset Y_{j}$$
for every $j\in\{1,\dots,i-1\}$.
\end{lemma}

\begin{proof}
There exists a geodesic $X_{0},\dots, X_{2i}$ in $\Gamma_{n-1}(\Pi)$ such
that $X_{2j}=Y_{j}$ for every $j\in\{0,\dots,i\}$.
Lemma \ref{lemma5-1red} gives the claim.
\end{proof}

\begin{lemma}\label{lemma5-3red}
If $m$ is an even integer not greater than $n$ then
for every isometric embedding $f$ of $\frac{1}{2}H_{m}$ in $\Gamma_{\delta}(\Pi)$
there is an $(n-m-1)$-dimensional singular subspace $M$ such that 
the image of $f$ is contained in $[M\rangle_{\delta}$,
in other words, $f$ can be considered as an isometric embedding of 
$\frac{1}{2}H_{m}$ in $\Gamma_{\delta}(\Pi_{M})$,
where $\Pi_{M}$ is the associated polar space of type {\rm $\textsf{D}_{m}$}.
\end{lemma}

\begin{proof}
Let $v$ and $w$ be opposite vertices of $\frac{1}{2}H_{m}$.
Then 
$$m=2d_{\delta}((f(v),f(w)))=d(f(v),f(w))=n-1-\dim(f(v)\cap f(w))$$
which implies that 
$$M:=f(v)\cap f(w)$$
is an $(n-m-1)$-dimensional singular subspace.
Our statement is a simple consequence of Lemmas \ref{lemma-half-cube-even} 
and \ref{lemma5-2red}.
\end{proof}

By Lemma \ref{lemma5-3red},
it is sufficient to prove Theorem \ref{theorem-apart-halfspin} only in the case when $m=n$.

\begin{rem}{\rm
In the case when  $m$ is odd, the graph $\frac{1}{2}H_{m}$ is more complicated
(Subsection 2.2) and we cannot prove the latter statement.
}\end{rem}

\subsection{Proof of Theorem \ref{theorem-apart-halfspin} for the case $m=n$}
Let $\Pi=(P,{\mathcal L})$ be a polar space of type $\textsf{D}_{n}$
and let ${\mathcal F}=\{p_{1},\dots,p_{2n}\}$ be a frame of $\Pi$.
For every integer $k\in \{1,\dots, n-1\}$ we denote by ${\mathcal A}_{k}$
the associated apartment of ${\mathcal G}_{k}(\Pi)$, i.e.
the set formed by all $k$-dimensional singular subspaces of type 
$$
\langle p_{i_{1}},\dots, p_{i_{k+1}}\rangle.
$$
Note that points $p_{i_{1}},\dots, p_{i_{k}}$ span a $(k-1)$-dimensional singular subspace 
if and only if 
$$\{i_{1},\dots,i_{k}\}\cap \{\sigma(i_{1}),\dots,\sigma(i_{k})\}=\emptyset,$$
or equivalently, these points are mutually collinear.
In this case, we will say that $i_{1},\dots,i_{k}$ form an {\it admissible} set.

Let $\delta\in \{+,-\}$.
Then 
$${\mathcal A}:={\mathcal A}_{n-1}\cap {\mathcal G}_{\delta}(\Pi)$$
is the apartment of ${\mathcal G}_{\delta}(\Pi)$ associated to ${\mathcal F}$.
For every admissible set $i_{1},\dots,i_{k}$ we denote by 
${\mathcal A}(i_{1},\dots,i_{k})$ the set of all elements of ${\mathcal A}$ containing 
$p_{i_{1}},\dots, p_{i_{k}}$, i.e.
$${\mathcal A}(i_{1},\dots,i_{k})={\mathcal A}\cap [S\rangle_{\delta},$$
where $S=\langle p_{i_{1}},\dots, p_{i_{k}}\rangle$.

It was noted in Subsection 4.2 that the restriction of $\Gamma_{\delta}(\Pi)$ to ${\mathcal A}$
is isomorphic to $\frac{1}{2}H_{n}$.
Let $f:{\mathcal A}\to {\mathcal G}_{\delta}(\Pi)$
be an isometric embedding of $\frac{1}{2}H_{n}$ in $\Gamma_{\delta}(\Pi)$.
Then 
$$\dim (X\cap Y)=\dim (f(X)\cap f(Y))$$
for all $X,Y\in {\mathcal A}$.

Using induction 
we show that ${\mathcal X}:=f({\mathcal A})$ is an apartment of ${\mathcal G}_{\delta}(\Pi)$ if $n$ is even. 
For $n=4$ this is true, see Example \ref{exmp3-1}. 
Suppose that the statement holds for $n=2l\ge 4$ and consider the case when $n=2(l+1)\ge 6$.

For every admissible set $i_{1},\dots,i_{k}$ we define
$${\mathcal X}(i_{1},\dots,i_{k}):=f({\mathcal A}(i_{1},\dots,i_{k})).$$
If $i,j$ is an admissible pair  then the restriction of $\Gamma_{\delta}(\Pi)$ to ${\mathcal A}(i,j)$
is isomorphic to $\frac{1}{2}H_{n-2}$ and
Lemma \ref{lemma5-3red} implies the existence of a line $L_{ij}$ such that 
$${\mathcal X}(i,j)\subset [L_{ij}\rangle_{\delta}.$$

\begin{lemma}\label{lemma5-4}
If $i,j$ and $s,t$ are admissible pairs satisfying
$\{i,j\}\cap \{s,t\}=\emptyset$
then 
$$L_{ij}\cap L_{st}=\emptyset.$$
\end{lemma}

\begin{proof}
We choose $X\in {\mathcal A}(i,j)$ and $Y\in {\mathcal A}(s,t)$
such that $X\cap Y=\emptyset$.
Then  
$$f(X)\cap f(Y)=\emptyset.$$
The statement follows from the fact that 
$f(X)\in {\mathcal X}(i,j)$ and $f(Y)\in {\mathcal X}(s,t)$ 
contain $L_{ij}$ and $L_{st}$, respectively. 
\end{proof}

Let $i,j$ be an admissible pair. 
Then  $\Pi_{p_{i}p_{j}}$ is a polar space of type $\textsf{D}_{n-2}$ 
and ${\mathcal A}(i,j)$ is an apartment of ${\mathcal G}_{\delta}(\Pi_{p_{i}p_{j}})$.
The inductive hypothesis implies that ${\mathcal X}(i,j)$ is an apartment of  ${\mathcal G}_{\delta}(\Pi_{L_{ij}})$.
Our embedding induces an isomorphism between 
the restrictions of $\Gamma_{\delta}(\Pi)$ to ${\mathcal A}(i,j)$ and ${\mathcal X}(i,j)$
(they are isomorphic to $\frac{1}{2}H_{n-2}$) and, 
by Subsection 4.2, we have the following possibilities:
\begin{enumerate}
\item[(1)] $f$ transfers every star of ${\mathcal A}(i,j)$ to a star of ${\mathcal X}(i,j)$,
\item[(2)] $n=6$ and $f$ sends every star of ${\mathcal A}(i,j)$ to a special subset of ${\mathcal X}(i,j)$.
\end{enumerate}

\begin{lemma}\label{lemma5-5}
The second possibility is not realized. 
\end{lemma}

\begin{proof}
Suppose that $n=6$ and $f$ transfers every star of ${\mathcal A}(i,j)$ to a special subset of ${\mathcal X}(i,j)$.
We take any $t$ such that $i,j,t$ form an admissible set.
Then ${\mathcal A}(i,j,t)$ is a star of ${\mathcal A}(i,j)$ and,
by our assumption, ${\mathcal X}(i,j,t)$ is a special subset of ${\mathcal X}(i,j)$.
The latter implies that the intersection of all elements from ${\mathcal X}(i,j,t)$
coincides with $L_{ij}$.

On the other hand, ${\mathcal A}(i,j,t)$ is a star in ${\mathcal A}(i,t)$ and ${\mathcal A}(j,t)$.
Hence ${\mathcal X}(i,j,t)$ is a special subset in ${\mathcal X}(i,t)$ and ${\mathcal X}(j,t)$.
By the above arguments, the lines $L_{it}$ and $L_{jt}$ both coincide with $L_{ij}$.

Now we take any $s\ne t$ such that $i,j,s$ form an admissible set.
Then ${\mathcal A}(i,j,s)$ is a star of ${\mathcal A}(i,j)$
and ${\mathcal X}(i,j,s)$ is a special subset of ${\mathcal X}(i,j)$.
As above, we establish that $L_{is}$ coincides with $L_{ij}$.
Thus $L_{is}=L_{jt}$ which contradict Lemma \ref{lemma5-4}.
\end{proof}

%As above, we suppose that $i,j$ is an admissible pair. 
Lemma \ref{lemma5-5} shows that
$f$ transfers every star of ${\mathcal A}(i,j)$ to a star of ${\mathcal X}(i,j)$.
Then, by Subsection 4.2,
for every $t$ such that $i,j,t$ is admissible set 
there is a plane $S_{ij}(t)$ containing $L_{ij}$ and such that 
${\mathcal X}(i,j,t)$ consists of all elements of ${\mathcal X}(i,j)$ containing $S_{ij}(t)$.
The intersection of all elements from ${\mathcal X}(i,j,t)$ coincides with $S_{ij}(t)$.
Since
$${\mathcal A}(i,j,t)=
{\mathcal A}(i,j)\cap {\mathcal A}(i,t)\cap {\mathcal A}(j,t),$$
we have
$${\mathcal X}(i,j,t)={\mathcal X}(i,j)\cap{\mathcal X}(i,t)\cap{\mathcal X}(j,t)$$
and
$$S_{ij}(t)=S_{it}(j)=S_{jt}(i).$$
In what follows this plane will be denoted by $S_{ijt}$.
The following fact is obvious.

\begin{lemma}\label{lemma5-6}
If $i,j$ is an admissible pair then all $S_{ijt}$ with $t\not\in \{i,j,\sigma(i),\sigma(j)\}$ 
form a frame of the polar space $\Pi_{L_{ij}}$
and ${\mathcal X}(i,j)$ is the apartment of ${\mathcal G}_{\delta}(\Pi_{L_{ij}})$ associated to this frame.
\end{lemma}

The plane $S_{ijt}$ contains the lines $L_{ij},L_{it},L_{jt}$. 
Hence for any pair of these lines one of the following possibilities is realized:
the lines are coincident or they are distinct and span $S_{ijt}$. 

\begin{lemma}\label{lemma5-7}
If $i,j,t$ is an admissible set 
then the lines $L_{ij},L_{it},L_{jt}$ are mutually distinct, i.e.
any two of these lines span the plane $S_{ijt}$.
\end{lemma}

\begin{proof}
Suppose that $L_{ij}=L_{it}$.
Consider $X\in {\mathcal A}(i,j)$ and $Y\in {\mathcal A}(i,t)$
such that 
$$X\cap Y=p_{i}p_{s}\;\mbox{ and }\;s\ne j,t.$$
Then $X,Y\in {\mathcal A}(i,s)$ and
$f(X),f(Y)\in {\mathcal X}(i,s)$. 
Since $f(X)\cap f(Y)$ is a line, we have
$$
f(X)\cap f(Y)=L_{is}.
$$
On the other hand, $f(X)\in {\mathcal X}(i,j)$ contains $L_{ij}$ and 
$f(Y)\in {\mathcal X}(i,t)$ contains $L_{it}$.
This means that the line $L_{ij}=L_{it}$ is contained in both $f(X)$ and $f(Y)$,
i.e. it coincides with $L_{is}$.
Therefore, the lines $L_{is}$ and $L_{jt}$ have a non-empty intersection
($L_{is}$ coincides with $L_{ij}$ and the latter line has a non-empty intersection with $L_{jt}$)
which is impossible by Lemma \ref{lemma5-4}.
\end{proof}

Consider the graph $\Gamma_{1}(\Pi)$ whose vertex set is ${\mathcal G}_{1}(\Pi)$
and two lines are adjacent vertices  of $\Gamma_{1}(\Pi)$
if they span a plane.
By \cite[Proposition 4.16]{Pankov1}, 
maximal cliques of this graph are of the following two types:
\begin{enumerate}
\item[$\bullet$] the {\it star} $[p, U]_{1}$, where $p$ is a point of $U\in {\mathcal G}_{n-1}(\Pi)$;
\item[$\bullet$] the {\it top} $\langle T]_{1}$, where $T\in {\mathcal G}_{2}(\Pi)$.
\end{enumerate}
An easy verification shows that the restriction of the graph $\Gamma_{1}(\Pi)$
to ${\mathcal A}_{1}$ has precisely the following two types of maximal cliques:
\begin{enumerate}
\item[$\bullet$] the {\it star} ${\mathcal A}_{1}\cap[p_{i}, U]_{1}$, where $U\in {\mathcal A}_{n-1}$ contains $p_{i}$;
\item[$\bullet$] the {\it top} ${\mathcal A}_{1}\cap \langle T]_{1}$, where $T\in {\mathcal A}_{2}$.
\end{enumerate}
Denote by $g$ the mapping of ${\mathcal A}_{1}$ to ${\mathcal G}_{1}(\Pi)$
sending every line $p_{i}p_{j}$ to the line $L_{ij}$.
It follows from Lemma \ref{lemma5-7} that
$g$ transfers adjacent  vertices of  $\Gamma_{1}(\Pi)$  to adjacent  vertices of  $\Gamma_{1}(\Pi)$.
Then $g$ sends every star of ${\mathcal A}_{1}$ to a subset in a star or a top of ${\mathcal G}_{1}(\Pi)$.
We cannot state that $g$ is injective, but its restriction to a star or a top of ${\mathcal A}_{1}$ is injective.

\begin{lemma}\label{lemma5-8}
Let $U\in {\mathcal A}$ and $p_{i}\in U$.
Then 
\begin{equation}\label{eq5-1}
g({\mathcal A}_{1}\cap [p_{i},U]_{1})
\end{equation}
is contained in a star of ${\mathcal G}_{1}(\Pi)$.
\end{lemma}

\begin{proof}
Let $J$ be the set of all $j\ne i$ such that $p_{j}\in U$. 
Then ${\mathcal A}_{1}\cap [p_{i},U]_{1}$ consists of all lines $p_{i}p_{j}$ with $j\in J$
and \eqref{eq5-1} is formed by all $L_{ij}$ with $j\in J$.
Now we fix $j\in J$ and consider $f(U)$ as an element of ${\mathcal X}(i,j)$.
Since $U$ belongs to ${\mathcal A}(i,j,t)$ for every $t\in J\setminus \{j\}$,
the planes $S_{ijt}$, $t\in J\setminus \{j\}$ are contained in $f(U)$  
and Lemma \ref{lemma5-6} implies that $f(U)$ is spanned by these planes.
Then Lemma \ref{lemma5-7} guarantees that 
$f(U)$ is spanned by all lines belonging to \eqref{eq5-1}. 
This means that \eqref{eq5-1} cannot be contained in a top 
and we get the claim.
\end{proof}

Let $U$ and $T$ be distinct elements of ${\mathcal A}$ containing $p_{i}$.
Lemma \ref{lemma5-8} implies the existence of points $q_{i}(U),q_{i}(T)$
and maximal singular subspaces $U',T'$ such that
$$g({\mathcal A}_{1}\cap [p_{i},U]_{1})\subset [q_{i}(U),U']_{1}$$
and
$$g({\mathcal A}_{1}\cap [p_{i},T]_{1})\subset [q_{i}(T),T']_{1}.$$
If $U$ and $T$ are adjacent vertices of the graph $\Gamma_{\delta}(\Pi)$
then 
$$\dim(U\cap T)=n-3\ge 3$$
and
$$({\mathcal A}_{1}\cap [p_{i},U]_{1})\cap({\mathcal A}_{1}\cap [p_{i},T]_{1})$$
contains more than one element.
The same holds for 
$$[q_{i}(U),U']_{1}\cap [q_{i}(T),T']_{1}$$
(since the restriction of $g$ to every star is injective).
This means that $q_{i}(U)=q_{i}(T)$
(if $q_{i}(U)$ and $q_{i}(T)$ are distinct then the latter intersection contains at most one element).
The graph $\Gamma_{\delta}(\Pi)$ is connected and
the equality $q_{i}(U)=q_{i}(T)$ holds for any $U,T\in {\mathcal A}$
containing $p_i$.

Since every line $p_{i}p_{j}$ is contained in a certain $U\in {\mathcal A}$,
we get the following:
for every point $p_{i}$ there is a point $q_{i}$ such that 
$$g({\mathcal A}_{1}\cap [p_{i}\rangle_{1})\subset [q_{i}\rangle_{1}.$$
The point $q_{i}$ is the intersection of all lines $L_{ij}$.

Show that $q_{i}\ne q_{j}$ if $i\ne j$.
The lines $L_{is}$ and $L_{jt}$ contain $q_{i}$ and $q_{j}$ (respectively)
and, by Lemma \ref{lemma5-4},
these lines are disjoint if $\{i,s\}\cap \{j,t\}=\emptyset$.

Thus $L_{ij}$ is spanned by $q_{i}$ and $q_{j}$.
Then for every $U\in {\mathcal A}$
the subspace $f(U)$ is spanned by all $q_{i}$ such that $p_{i}\in U$
(see the proof of Lemma \ref{lemma5-8}). 
So, the required statement is a simple consequence of the following.

\begin{lemma}\label{lemma5-9}
The points $q_{1},\dots,q_{2n}$ form a frame of $\Pi$.
\end{lemma}

To prove Lemma \ref{lemma5-9}
we need some elementary properties of collinearity relation in polar spaces.

\begin{lemma}\label{lemma5-10}
The following assertions are fulfilled:
\begin{enumerate}
\item[(1)] if $p\in P$ and $X\subset P$ then $p\perp X$ implies that $p\perp \langle X\rangle$,
\item[(2)] if $p\in P$ and $S$ is a maximal singular subspace then
$p\perp S$ implies that $p\in S$.
\end{enumerate}
\end{lemma}
\begin{proof}
See, for example, \cite[Section 4.1]{Pankov1}.
\end{proof}

\begin{proof}[Proof of Lemma \ref{lemma5-9}]
We have $q_{i}\perp q_{j}$ if $j\ne \sigma(i)$ (since these points span the line $L_{ij}$).
Suppose that $q_{i}\perp q_{\sigma(i)}$ for a certain $i$.
Then $q_{\sigma(i)}$ is collinear to all points $q_{j}$.
Every element of ${\mathcal X}$ is spanned by some $q_{j}$
and the first statement of Lemma \ref{lemma5-10}  guarantees that
$q_{\sigma(i)}\perp S$ for every $S\in {\mathcal X}$.
By the second statement of Lemma \ref{lemma5-10},  
 every element of ${\mathcal X}$ contains $q_{\sigma(i)}$ and
 any two elements of ${\mathcal X}$ have a non-empty intersection.
The latter contradicts the fact that ${\mathcal X}$ is 
the image of an isometric embedding of $\frac{1}{2}H_{n}$  in $\Gamma_{\delta}(\Pi)$.
\end{proof}

\section{Conjecture}
In this section we will suppose that 
$\Pi=(P,{\mathcal L})$ is a polar space of type $\textsf{D}_n$, 
where $n$ is an odd integer greater than $4$.
Let ${\mathcal A}$ be the apartment of ${\mathcal G}_{\delta}(\Pi)$, $\delta\in \{+,-\}$
defined by a frame ${\mathcal F}=\{p_{1},\dots,p_{2n}\}$.
Let also $f:{\mathcal A}\to {\mathcal G}_{\delta}(\Pi)$ be 
an isometric embedding of $\frac{1}{2}H_{n}$ in $\Gamma_{\delta}(\Pi)$.
For every $i\in \{1,\dots,2n\}$ the restriction of $\Gamma_{\delta}(\Pi)$ to 
${\mathcal A}\cap [p_{i}\rangle_{\delta}$ is
isomorphic to $\frac{1}{2}H_{n-1}$ and Lemma \ref{lemma5-3red} implies the existence 
of a point $q_{i}$ such that 
$$f({\mathcal A}\cap [p_{i}\rangle_{\delta})\subset [q_{i}\rangle_{\delta}.$$
We can prove that there are precisely the following possibilities:
\begin{enumerate}
\item[(1)] all $q_{i}$ are coincident,
\item[(2)]  all $q_{i}$ are mutually distinct.
\end{enumerate}
In the second case, an easy verification shows that all $q_{i}$ form a frame 
and $f({\mathcal A})$ is the apartment of ${\mathcal G}_{\delta}(\Pi)$
associated to this frame.

We conjecture that the first possibility is realized, i.e. 
there exists an isometric embedding of $\frac{1}{2}H_{n}$ in $\Gamma_{\delta}(\Pi)$
whose image is not an apartment of  ${\mathcal G}_{\delta}(\Pi)$.

In the first case,  $f$ can be considered as an isometric embedding of $\frac{1}{2}H_{n}$
in $\Gamma_{\delta}(\Pi_{q})$, where $q=q_{1}=\dots=q_{2n}$.
The rank of the polar space $\Pi_{q}$ is equal to $n-1$.
Thus the graphs $\frac{1}{2}H_{n}$ and $\Gamma_{\delta}(\Pi_{q})$ are of the same diameter
$\frac{n-1}{2}$. This trivial observation supports our conjecture.

\section{Isometric embeddings of half-spin Grassmann graphs}

\subsection{Application of Theorem \ref{theorem-apart-halfspin}}
Let $\Pi=(P,{\mathcal L})$ and $\Pi'=(P',{\mathcal L}')$ be polar spaces of rank $n$ and $n'$,
respectively.
The existence of isometric embeddings of $\Gamma_{n-1}(\Pi)$ in $\Gamma_{n'-1}(\Pi')$
implies that the diameter of $\Gamma_{n-1}(\Pi)$ is not greater than the diameter of $\Gamma_{n'-1}(\Pi')$,
i.e. $n\le n'$.

First, we consider the case when $n=n'$.
Let $f:P\to P'$ be a mapping which transfers every frame of $\Pi$ to a frame of $\Pi'$.
By \cite[Subsection 4.9.6]{Pankov1}, $f$ sends lines of $\Pi$ to subsets contained in lines of $\Pi'$;
moreover, if $S$ is a singular subspace of $\Pi$
then $f(S)$ spans a singular subspace  of $\Pi'$ whose dimension is equal to the dimension of $S$.
The mapping which transfers every $X\in {\mathcal G}_{n-1}(\Pi)$ to $\langle f(X)\rangle$
is an isometric embedding of $\Gamma_{n-1}(\Pi)$ in $\Gamma_{n-1}(\Pi')$.
If $\Pi$ and $\Pi'$ both are polar spaces of type $\textsf{D}_{n}$ then 
the restriction of this mapping to ${\mathcal G}_{\delta}(\Pi)$, $\delta \in \{+,-\}$
is an isometric embedding of $\Gamma_{\delta}(\Pi)$ in $\Gamma_{\gamma}(\Pi')$,
$\gamma \in \{+,-\}$.

Suppose that $n\le n'$ and $M$ is an $(n'-n-1)$-dimensional singular subspace of $\Pi'$.
Then $\Pi'_{M}$ is a polar space of rank $n$ 
(if $n=n'$ then $M=\emptyset$ and $\Pi'_{M}$ coincides with $\Pi'$). 
Every frame preserving mapping of $\Pi$ to $\Pi'_{M}$ induces an 
isometric embedding of $\Gamma_{n-1}(\Pi)$ in $\Gamma_{n-1}(\Pi'_{M})$.
This is an isometric embedding of $\Gamma_{n-1}(\Pi)$ in $\Gamma_{n'-1}(\Pi')$,
since $\Gamma_{n-1}(\Pi'_{M})$ can be naturally identified with 
the restriction of $\Gamma_{n'-1}(\Pi')$ to $[M\rangle_{n'-1}$.
In the case when $\Pi$ and $\Pi'$ are polar spaces of type $\textsf{D}_{n}$ and $\textsf{D}_{n'}$
(respectively), we get an isometric embedding of 
$\Gamma_{\delta}(\Pi)$, $\delta \in \{+,-\}$ in $\Gamma_{\gamma}(\Pi')$, $\gamma \in \{+,-\}$.

\begin{theorem}[\cite{Pankov3}]\label{theorem-dual-polar-emb}
Every isometric embedding of $\Gamma_{n-1}(\Pi)$ in $\Gamma_{n'-1}(\Pi')$
is induced by a frame preserving mapping of $\Pi$ to $\Pi'_{M}$,
where $M$ is an $(n'-n-1)$-dimensional singular subspace of $\Pi'$.
\end{theorem}

As an application of  Theorem \ref{theorem-apart-halfspin} we prove the following.

\begin{theorem}\label{theorem-half-spin-emb}
Suppose that $\Pi$ and $\Pi'$ are polar spaces of type {\rm $\textsf{D}_{n}$} and 
{\rm $\textsf{D}_{n'}$}, respectively.
If $n$ is an even integer not less than $6$ then every isometric embedding of 
$\Gamma_{\delta}(\Pi)$, $\delta \in \{+,-\}$ in $\Gamma_{\gamma}(\Pi')$, $\gamma \in \{+,-\}$
is induced by a frame preserving mapping of $\Pi$ in $\Pi'_{M}$, where 
$M$ is an $(n'-n-1)$-dimensional singular subspace of $\Pi'$.
\end{theorem}

\subsection{Proof of Theorem \ref{theorem-half-spin-emb}}
Suppose that $f:{\mathcal G}_{\delta}(\Pi)\to {\mathcal G}_{\gamma}(\Pi')$ is an isometric embedding of 
$\Gamma_{\delta}(\Pi)$ in $\Gamma_{\gamma}(\Pi')$ and $n$ is even.
Then $f$ transfers maximal cliques of $\Gamma_{\delta}(\Pi)$ (stars and special subspaces) 
to subsets in maximal cliques of $\Gamma_{\gamma}(\Pi')$.

\begin{lemma}\label{lemma7-1}
The following assertions are fulfilled:
\begin{enumerate}
\item[{\rm (1)}] distinct maximal cliques go to subsets of distinct maximal cliques,
\item[{\rm (2)}] the image of every maximal clique of $\Gamma_{\delta}(\Pi)$ is contained in unique 
maximal clique of $\Gamma_{\gamma}(\Pi')$.
\end{enumerate}
\end{lemma}

\begin{proof}
(1). If ${\mathcal X}$ and ${\mathcal Y}$ are distinct maximal cliques of $\Gamma_{\delta}(\Pi)$
whose images are contained in the same maximal clique of $\Gamma_{\gamma}(\Pi')$
then there exist non-adjacent vertices $X\in {\mathcal X}$ and $Y\in {\mathcal Y}$
such that $f(X)$ and $f(Y)$ are adjacent vertices of $\Gamma_{\gamma}(\Pi')$.
This contradicts the fact that $f$ is an embedding.

(2).
Let ${\mathcal C}$ be a maximal clique of $\Gamma_{\delta}(\Pi)$
and let $S$ be the associated singular subspace of $\Pi$
(an element of ${\mathcal G}_{n-4}(\Pi)$ or ${\mathcal G}_{-\delta}(\Pi)$).
By \cite[Proposition 4.7]{Pankov1}, there is a frame ${\mathcal F}$ 
such that $S$ is spanned by a subset of ${\mathcal F}$.
Let ${\mathcal A}$ be the apartment of ${\mathcal G}_{\delta}(\Pi)$
defined by ${\mathcal F}$.
Then ${\mathcal A}\cap {\mathcal C}$ is a base of the projective space ${\mathcal C}$
(see Subsection 4.2).
Theorem \ref{theorem-apart-halfspin} states that 
$f({\mathcal A})$ is an apartment of ${\mathcal G}_{\delta}(\Pi'_{M})$,
where $M$ is an $(n'-n-1)$-dimensional singular subspace of $\Pi'$.
This guarantees that $f({\mathcal A}\cap {\mathcal C})$
is an independent subset in every maximal clique of $\Gamma_{\gamma}(\Pi')$
containing $f({\mathcal C})$
and  the dimension of the subspace spanned by $f({\mathcal A}\cap {\mathcal C})$ is not less than $3$.

Now suppose that $f({\mathcal C})$ is contained in two distinct maximal cliques of $\Gamma_{\gamma}(\Pi')$.
If these cliques are of the same type 
then $f({\mathcal C})$ is contained in a line of  ${\mathcal G}_{\gamma}(\Pi')$. 
If one of the cliques is a star and the other is a special subspace then
$f({\mathcal C})$ spans a proper subspace in the star and the dimension of this subspace is not greater
than $2$.  Each of these possibilities contradicts the fact established above.
\end{proof}

\begin{lemma}\label{lemma7-2}
One of the following possibilities is realized:
\begin{enumerate}
\item[{\rm (A)}] 
stars go to subsets of stars and special subspaces go to subsets of special subspaces, 
\item[{\rm (B)}]
stars go to subsets of special subspaces and special subspaces go to subsets of stars. 
\end{enumerate}
The second possibility can be realized only in the case when $n=4$.
\end{lemma}

\begin{proof}
Let ${\mathcal C}$ and ${\mathcal C}'$ be  maximal cliques of $\Gamma_{\delta}(\Pi)$
such that ${\mathcal C}\cap {\mathcal C}'$ is a plane.
Then one of these cliques is a star $[M\rangle_{\delta}$, $M\in {\mathcal G}_{n-4}(\Pi)$
and the other is a special subspace $[U]_{\delta}$, $U\in {\mathcal G}_{-\delta}(\Pi)$.
There is a frame ${\mathcal F}$ 
such that $M$ and $U$ are spanned by subsets of ${\mathcal F}$ \cite[Proposition 4.7]{Pankov1}.
Let ${\mathcal A}$ be the apartment of ${\mathcal G}_{\delta}(\Pi)$ defined by ${\mathcal F}$.
Then $${\mathcal A}\cap {\mathcal C}\cap {\mathcal C}'$$ is a base of the plane ${\mathcal C}\cap {\mathcal C}'$.
By Theorem \ref{theorem-apart-halfspin},
$f({\mathcal A})$ is an apartment of ${\mathcal G}_{\delta}(\Pi'_{M})$,
where $M$ is an $(n'-n-1)$-dimensional singular subspace of $\Pi'$. 
This means that 
$$f({\mathcal A}\cap {\mathcal C}\cap {\mathcal C}')$$
is an independent subset in the intersection of the maximal cliques of $\Gamma_{\gamma}(\Pi')$ 
containing  $f({\mathcal C})$ and $f({\mathcal C}')$. 
Since this set contains three elements, the dimension of the intersection is not less than $2$. 
Therefore, the intersection of the maximal cliques of $\Gamma_{\gamma}(\Pi')$ 
containing  $f({\mathcal C})$ and $f({\mathcal C}')$ is a plane.

Lemma \ref{lemma7-1} implies that $f$ induces an injective mapping 
of the set of maximal cliques of $\Gamma_{\delta}(\Pi)$ to 
the set of maximal cliques of $\Gamma_{\gamma}(\Pi')$.
Using Lemma \ref{lemma4-1}, we show that one of the possibilities (A) or (B) is realized.

Now suppose that $n\ge 6$.
Let $U\in {\mathcal G}_{-\delta}(\Pi)$. We take a frame ${\mathcal F}$
such that $U$ is spanned by a subset of ${\mathcal F}$
and denote by ${\mathcal A}$ the associated apartment of ${\mathcal G}_{\delta}(\Pi)$.
Then ${\mathcal A}\cap [U]_{\delta}$ is a base of the projective space $[U]_{\delta}$.
As above, we establish that 
$f({\mathcal A}\cap [U]_{\delta})$
is an independent subset in the maximal clique of $\Gamma_{\gamma}(\Pi')$
containing $f([U]_{\delta})$. 
Since $f({\mathcal A}\cap [U]_{\delta})$ consists of $n$ elements and $n>4$,
the maximal clique containing $f([U]_{\delta})$ cannot be a star
and the possibility (B) is not realized.
 \end{proof}

 From this moment we suppose that $f$ satisfies (A).
Then for every $U\in {\mathcal G}_{-\delta}(\Pi)$ there is unique $U'\in {\mathcal G}_{-\gamma}(\Pi')$
such that 
$$f([U]_{\delta})\subset [U']_{\gamma}.$$ 
We set $f(U):=U'$ 
%for every $U\in {\mathcal G}_{-\delta}(\Pi)$ 
and extend $f$ to a mapping of ${\mathcal G}_{n-1}(\Pi)$ to ${\mathcal G}_{n'-1}(\Pi')$.
The first part of Lemma \ref{lemma7-1} guarantees that  
the extension  is injective.  

\begin{lemma}\label{lemma7-3}
The following assertions are fulfilled:
\begin{enumerate}
\item[{\rm (1)}] $f$ transfers adjacent vertices of $\Gamma_{n-1}(\Pi)$ 
to adjacent vertices of $\Gamma_{n'-1}(\Pi')$,
\item[{\rm (2)}] $d(X,Y)=d(f(X),f(Y))$ if $X,Y\in {\mathcal G}_{\delta}(\Pi)$.
\end{enumerate}
 \end{lemma}

\begin{proof}
(1).
If $X\in {\mathcal G}_{\delta}(\Pi)$ and $Y\in {\mathcal G}_{-\delta}(\Pi)$
are adjacent vertices of $\Gamma_{n-1}(\Pi)$ then 
$X$ belongs to $[Y]_{\delta}$ and 
$$f(X)\in [f(Y)]_{\gamma}$$ 
which implies that $f(X)$ and $f(Y)$ are adjacent vertices of $\Gamma_{n'-1}(\Pi')$.

(2). This follows from the fact that 
the restriction of $f$ to ${\mathcal G}_{\delta}(\Pi)$ is 
an isometric embedding of $\Gamma_{\delta}(\Pi)$ in $\Gamma_{\gamma}(\Pi')$.
\end{proof}

By Theorem \ref{theorem-dual-polar-emb}, it is sufficient to show that $f$
is an isometric embedding of $\Gamma_{n-1}(\Pi)$ in $\Gamma_{n'-1}(\Pi')$.
Since for any two distinct elements of ${\mathcal G}_{n-1}(\Pi)$ there is an apartment containing them
\cite[Proposition 4.7]{Pankov1},
we need to show that the restriction of $f$ to every apartment of ${\mathcal G}_{n-1}(\Pi)$
is an isometric embedding of $H_{n}$ in $\Gamma_{n'-1}(\Pi')$.

Let ${\mathcal A}$ be an apartment of ${\mathcal G}_{n-1}(\Pi)$.
Then 
$${\mathcal A}_{\delta}:={\mathcal A}\cap {\mathcal G}_{\delta}(\Pi)$$ 
is an apartment of ${\mathcal G}_{\delta}(\Pi)$
and ${\mathcal A}\setminus {\mathcal A}_{\delta}$
is an apartment of ${\mathcal G}_{-\delta}(\Pi)$.
Every special subset of ${\mathcal A}_{\delta}$ is defined by an element of 
${\mathcal A}\setminus {\mathcal A}_{\delta}$,
i.e. it is the intersection of ${\mathcal A}_{\delta}$ and $[U]_{\delta}$,
where $U\in{\mathcal A}\setminus {\mathcal A}_{\delta}$. 
Theorem \ref{theorem-apart-halfspin} states that 
$f({\mathcal A}_{\delta})$ is an apartment of ${\mathcal G}_{\gamma}(\Pi'_{M})$,
where $M$ is an $(n'-n-1)$-dimensional singular subspace of $\Pi'$.
Let ${\mathcal A}'$ be the apartment of ${\mathcal G}_{n-1}(\Pi'_{M})$ containing $f({\mathcal A}_{\delta})$.
Then ${\mathcal A}'\setminus f({\mathcal A}_{\delta})$ is an apartment of 
${\mathcal G}_{-\gamma}(\Pi'_{M})$ and every special subset of $f({\mathcal A}_{\delta})$
is defined by an element of ${\mathcal A}'\setminus f({\mathcal A}_{\delta})$.
Since $f$ transfers special subsets of ${\mathcal A}_{\delta}$ to special subsets of  $f({\mathcal A}_{\delta})$,
we have 
$$f({\mathcal A}\setminus {\mathcal A}_{\delta})={\mathcal A}'\setminus f({\mathcal A}_{\delta})$$
which implies that $f({\mathcal A})={\mathcal A}'$.

Now we show that 
\begin{equation}\label{eq7-1}
d(X,Y)=d(f(X),f(Y))
\end{equation}
for all $X,Y\in {\mathcal A}$.
By the second part of Lemma \ref{lemma7-3}, this is true if $X,Y\in {\mathcal A}_{\delta}$.
In the case when $X\in {\mathcal A}\setminus{\mathcal A}_{\delta}$ and $Y\in {\mathcal A}_{\delta}$,
we have $d(X,Y)\le n-1$ and
there is a geodesic 
$Y_{0},Y_{1},\dots, Y_{i}$ in $\Gamma_{n-1}(\Pi)$ such that 
$$Y_{0}\in {\mathcal A}_{\delta},\;\;Y_{1}=X,\;\;Y_{i}=Y.$$
Lemma \ref{lemma7-3} shows that
$$d(Y_{0},Y_{i})=d(f(Y_{0}),f(Y_{i}))$$
and
$$f(Y_{0}),f(Y_{1}),\dots, f(Y_{i})$$
is a geodesic in $\Gamma_{n'-1}(\Pi')$ which implies \eqref{eq7-1}.
If $X,Y\in{\mathcal A}\setminus{\mathcal A}_{\delta}$ and $d(X,Y)\le n-2$ then 
there is a geodesic
$Y_{0},Y_{1},\dots, Y_{i}$ in $\Gamma_{n-1}(\Pi)$ such that
$$Y_{0},Y_{i}\in {\mathcal A}_{\delta},\;\;Y_{1}=X\;\;Y_{i-1}=Y.$$
As above, the image of this geodesic is a geodesic in $\Gamma_{n'-1}(\Pi)$
and we get \eqref{eq7-1}.

Consider the case when $X,Y\in{\mathcal A}\setminus{\mathcal A}_{\delta}$ are opposite vertices of 
$\Gamma_{n-1}(\Pi)$.
It is clear that 
$${\mathcal A}\setminus\{Y\}=\{Z\in {\mathcal A}:d(X,Z)\le n-1\}$$ 
There is unique $Y'\in {\mathcal A}'$ satisfying $d(f(X),Y')=n$
and
$${\mathcal A}'\setminus \{Y'\}=\{Z'\in {\mathcal A}':d(f(X),Z')\le n-1\}.$$ 
It follows from the above arguments that
$$f({\mathcal A}\setminus\{Y\})={\mathcal A}'\setminus \{Y'\}.$$
Then $f(Y)=Y'$ (since $f({\mathcal A})={\mathcal A}'$)
and \eqref{eq7-1} holds.

So, we get the following.

\begin{prop}\label{prop7-1}
If $\Pi$ and $\Pi'$ are polar spaces of type {\rm $\textsf{D}_{n}$} and 
{\rm $\textsf{D}_{n'}$}{\rm (}respectively{\rm )}
and $n$ is even then every isometric embedding of 
$\Gamma_{\delta}(\Pi)$, $\delta \in \{+,-\}$ in $\Gamma_{\gamma}(\Pi')$, $\gamma \in \{+,-\}$
satisfying {\rm (A)}
is induced by a frame preserving mapping of $\Pi$ in $\Pi'_{M}$, where 
$M$ is an $(n'-n-1)$-dimensional singular subspace of $\Pi'$.
\end{prop}

Theorem \ref{theorem-half-spin-emb} is a direct consequence of 
Proposition \ref{prop7-1} and Lemma \ref{lemma7-2}.

\subsection{Remarks on the case $n=4$}
Suppose that $\Pi$ is a polar space of type $\textsf{D}_4$.
Then the associated half-spin Grassmannians (together with the families of lines)
also are polar spaces of type $\textsf{D}_4$.
These polar spaces are denoted by $\Pi_{+}$ and $\Pi_{-}$ (Example \ref{exmp3-1}).

Every maximal singular subspace of the polar space $\Pi_{\delta}$, $\delta\in \{+,-\}$
is a star (defined by a point of $\Pi$) or a special subspace
(associated to an element of ${\mathcal G}_{-\delta}(\Pi)$).
Using the intersection properties given in Section 4.1,
we can show that one of the half-spin Grassmannians of $\Pi_{\delta}$ consists of all stars 
and the other is formed by all special subspaces.
The corresponding polar spaces can be identified with $\Pi$ and $\Pi_{-\delta}$.
By Remark \ref{rem-halfspin}, they are isomorphic.

 So, the polar spaces $\Pi$ and $\Pi_{\delta}$, $\delta\in \{+,-\}$ are isomorphic.
Every collineation between $\Pi$ and $\Pi_{-\delta}$ induces
a collineation of $\Pi_{\delta}$ to itself (an automorphism of $\Gamma_{\delta}(\Pi)$)
which transfers stars to special subspaces and special subspaces to stars.

Let $f$ be an isometric embedding of 
$\Gamma_{\delta}(\Pi)$, $\delta \in \{+,-\}$ in $\Gamma_{\gamma}(\Pi')$, $\gamma \in \{+,-\}$
satisfying {\rm (B)}. 
We take any automorphism $h$ of $\Gamma_{\delta}(\Pi)$
sending stars to special subspaces and special subspaces to stars. 
Then $fh$ is an isometric embedding of 
$\Gamma_{\delta}(\Pi)$ in $\Gamma_{\gamma}(\Pi')$
satisfying {\rm (A)} and we can apply Proposition \ref{prop7-1}.

\end{document}